\documentclass[11pt]{article}
\usepackage{amsfonts,amssymb,amsmath,a4wide}
\usepackage[francais]{babel}
\def\Zm{{\mathbb Z}}

\def\Rm{{\mathbb R}}
\def\Tm{{\mathbb T}}
\def\Nm{{\mathbb N}}

\def\Qm{{\mathbb Q}}
\def\lto{\longrightarrow}
\def\lmto{\longmapsto}
\def\eq{\Longleftrightarrow}
\def\leq{\leqslant}
\def\geq{\geqslant}

\makeatletter
\renewcommand{\section}{\@startsection
{section}
{1}
{0mm}
{-1.2\baselineskip}
{\baselineskip}
{\center \scshape}}

\renewcommand{\subsection}{\@startsection
{subsection}
{2}
{0mm}
{-\baselineskip}
{0mm}
{\normalfont \normalsize \bfseries}}
\makeatother

\author{Patrick Bernard }
\title{Un r\'esultat de transfert en approximation
diophantienne}
\begin{document}

\begin{center}
\Large
\begin{scshape}
Une propri\'et\'e de transfert en approximation diophantienne\\
\vspace{.5cm}
\end{scshape}
\normalsize
Patrick BERNARD
\footnote{
Patrick Bernard, Institut Fourier, BP 74, 38402 Saint Martin
d'Hy\`eres cedex, FRANCE.
\\ Patrick.Bernard@ujf-grenoble.fr}\\
\vspace{.5cm}
janvier 2001, corrig\'e 
\'et\'e 2003
\end{center}

\vspace{.5cm}

\textit{
\textbf{R\'esum\'e : }
\'Etant donn\'e un vecteur $\omega \in \Rm^n$,
on d\'efinit la suite $T_i$ des p\'eriodes de $\omega$ comme la suite
des temps de meilleur retour pr\`es de l'origine
de la translation $x\lmto x+\omega$ sur le tore $\Tm ^n$.
On \'etudie comment les propri\'et\'es diophantiennes 
du vecteur $\omega$ peuvent \^etre exprim\'ees 
\`a l'aide de la suite des p\'eriodes. 
Plus pr\'ecis\'ement, on montre que 
si le vecteur $\omega$ est non r\'esonant, et si ses
p\'eriodes v\'erifient l'in\'egalit\'e 
$T_{i+1} \leq CT_i^{1+\tau}$ avec
$\tau<(n-1)^{-1}$, alors le 
vecteur $\omega$ est diophantien.}

\vspace{.5cm}

\textit{
\textbf{Abstract :}
Given a vector  $\omega \in \Rm^n$,
the sequence  $T_i$ of periods is defined as the sequence of times
of best returns near the origin of the translation
 $x\lmto x+\omega$ on the torus  $\Tm ^n$.
In the present paper, we study how the Diophantine
properties of $\omega$ can be expressed considering the
sequence of its periods.
More precidely, we prove that,
if the vector $\omega$ is not resonant,
and if the sequence of periods satisfy the inequality
$T_{i+1} \leq CT_i^{1+\tau}$ with
$\tau<(n-1)^{-1}$, then the 
vector $\omega$ is Diophantine.}

\vspace{.5cm}

\section*{Introduction}

\subsection{ }
Soit $\omega$ un vecteur de $\Rm^n$.
On note $\Tm^n$ le tore $(\Rm/\Zm)^n$, et 
$
t_{\omega}:\Tm^n \lto \Tm^n
$
le diff\'eomor\-phisme  $x \lmto x+\omega$.
On  s'int\'eresse au  syst\`eme dynamique discret 
engendr\'e par $t_{\omega}$.
On dit que $\omega$ est r\'esonant si il existe $k\in \Zm^n$
tel que 
le produit scalaire $\langle \omega, k \rangle $
est entier $(\in \Zm)$. Dans ce cas, 
le tore $\Tm^n$ est fibr\'e en tores de plus petite dimension
invariants par $t_{\omega}$.
Sinon, on dit que $\omega$ est non r\'esonant, et 
$t_{\omega}$ est uniquement ergodique et minimal,
c'est \`a dire que $\Tm^n$ est le seul compact non vide invariant
par $t_{\omega}$, et que la mesure de Haar est la seule 
mesure de probabilit\'e invariante.

Il est tr\`es utile en th\'eorie de la moyennisation 
(th\'eorie KAM par exemple)
de quantifier le caract\`ere non r\'esonant
de $\omega$.
Le plus simple 
consiste \`a minorer les ''petits diviseurs''
$d(\langle k,\omega\rangle,\Zm)$,  $k\in\Zm^n$.
On parle alors de propri\'et\'es diophantiennes lin\'eaires.
Cette approche est
assez naturelle, car les petits diviseurs apparaissent
directement au d\'enominateur des coefficients des s\'eries
de perturbation dans de nombreux probl\`emes.
Intuitivement, ceci revient a quantifier l'ergodicit\'e de
$t_{\omega}$, mais le sens dynamique des petits diviseurs n'est
pas tr\`es clair.

On peut aussi mesurer la
qualit\'e de l'approximation de 
$\omega$ par des vecteurs rationnels, c'est \`a dire estimer
$d(T\omega,\Zm^n)$, $T\in\Zm$.
On parle alors de propri\'et\'es diophantiennes simultan\'ees.
Cette approche n'a \'et\'e utilis\'ee que beaucoup plus 
r\'ecemment dans les probl\`emes de moyennisation,
voir \cite{Lochak}, o\`u elle permet une simplification
remarquable du th\'eor\`eme de Nekhoroshev.
On peut enfin se restreindre \`a la suite des meilleures 
approximations de $\omega$,
comme sugg\'er\'e dans \cite{Lochak}.
On d\'efinit pour ceci  les p\'eriodes $T_i(\omega)$
de  $\omega\in\Rm^n$
en posant $T_0(\omega)=1$
et 
\begin{equation}\label{periodes} 
T_{i+1}(\omega)=\min \big\{ T\in \Nm \text{ tel que } 
\|T \omega \|_{\Zm}< \|T_i(\omega) \omega \|_{\Zm} \big \}.
\end{equation}
Voir \ref{notations} ci-dessous pour les notations.
La croissance des p\'eriodes d\'etermine de mani\`ere tr\`es
condens\'ee les propri\'et\'es diophantiennes de $\omega$.

Un int\'er\^et majeur de l'approximation simultan\'ee est 
son caract\`ere profond\'ement dynamique.
En effet, si l'on consid\`ere un hom\'eomorphisme $\phi$
d'un espace m\'etrique et une orbite r\'ecurrente 
$\phi^n(x)$ de cet hom\'eomorphisme,
on peut d\'efinir les 
\'ecarts $d(\phi^n(x),x)$ et les p\'eriodes 
$T_i$ de cette orbite par $T_0=1$ et
$$
T_{i+1}=\min \big\{ T\in \Nm \text{ tel que } 
 d(\phi^{T}(x),x)<d(\phi^{T_i}(x),x)\big \}.
$$
Ces d\'efinition co\" \i ncident avec 
les notions vues au dessus
lorsque $\phi$ est la translation sur le tore de vecteur $\omega$.

La correspondance pr\'ecise entre les diff\'erents types de
propri\'et\'es  diophantiennes n'est pas imm\'ediate,
bien que ces propri\'et\'es soient  visiblement de m\^eme 
nature. 
Nous rappelons ici le tr\`es classique th\'eor\`eme de transfert
de Khintchine, et l'\'etendons aux propri\'et\'es diophantiennes
d\'efinies en terme de croissance des p\'eriodes,
qui n'avaient pas \'et\'e \'etudi\'ees jusqu'ici.
Cet article doit beaucoup \`a mes conversations avec Pierre 
Lochak
sur le lien entre 
les m\'ethodes de perturbation par approximation
simultan\'ees d\'evelopp\'ees dans \cite{Lochak} et 
la th\'eorie KAM. 
Je tiens aussi a remercier Xavier Buff,
qui m'a sugg\'er\'e une am\'elioration de la Propri\'et\'e 1.3,
ainsi que le \textit{referee} pour sa lecture attentive
qui a permis de nombreuses am\'eliorations du texte.

\subsection { }\label{notations}
\textsc{Notations :}
On utilise pour l'essentiel les notations de l'appendice 1
de \cite{Lochak}.
Soit $\omega=(\omega_1,\ldots,\omega_n)\in\Rm^n$,
on note
$$
 |\omega| 
 =\sup_i |\omega_i|.
$$
Si $\omega'=(\omega'_1,\ldots,\omega'_n)\in\Rm^n$,
on note
$$
\langle \omega,\omega'\rangle
=\omega_1\omega'_1+\ldots +\omega_n\omega'_n
$$
le produit scalaire standard,
et
$$\|\omega \|=
\sqrt{\langle \omega,\omega\rangle}.$$
Pour $\omega\in \Rm^n$, on note 
$$\| \omega\|_{\Zm} =
\min_{k\in \Zm^n} |\omega-k|,
$$
et de la m\^eme fa\c con, 
$\|x\|_{\Zm} =
\min_{k\in \Zm} |x-k| 
$
pour $x$ r\'eel.

\section{Th\'eor\`emes de Dirichlet}\label{Dirichlet}

Il est utile pour fixer les id\'ees de rappeler les deux
r\'esultats les plus simples de l'approximation diophantienne.
On pourra consulter \cite{Schmidt} pour les preuves.
\subsection{ }
\textsc{Th\'eor\`eme de Dirichlet : }\label{sim}
Consid\'erons un vecteur $\omega\in \Rm^n$, 
pour tout r\'eel $Q>0$, 
il existe un entier positif $T<Q$  
tel que
$$
\|T\omega\|_{\Zm} \leq Q^{-\frac{1}{n}}.
$$
En cons\'equence,
il existe une infinit\'e d'entiers $T$ tels que 
$$
\|T\omega\|_{\Zm} <  T^{-\frac{1}{n}}.
$$

\subsection{ }
\textsc{Th\'eor\`eme : }\label{lin}
Consid\'erons un vecteur $\omega\in \Rm^n$, 
pour tout r\'eel $Q>0$
il existe un vecteur $k\in\Zm^n$
v\'erifiant $|k|<Q$
et tel que
$$
\| \langle k, \omega \rangle\|_{\Zm}
\leq
Q^{-n}.
$$
En cons\'equence, 
il existe une infinit\'e de vecteurs entiers $k$ tels que 
$$
\| \langle k, \omega \rangle\|_{\Zm}
<
|k |^{-n}.
$$%
\subsection{ }
\textsc{Propri\'et\'e : }
La suite $T_i(\omega)$ des p\'eriodes de $\omega$
satisfait:
$$
\frac{1}{T_i(\omega)+T_{i+1}(\omega)}
\leq
\|T_i(\omega)\omega\|_{\Zm}
\leq
\frac{1}{T_{i+1}(\omega)^{1/n}}.
$$
L'in\'egalit\'e de droite est une cons\'equence 
directe du th\'eor\`eme
de Dirichlet avec $Q=T_{i+1}(\omega)$.
En effet, on a alors 
$$
\|T_i(\omega)\omega\|_{\Zm}=
\min _{1\leq T<T_{i+1}(\omega)}\|T\omega\|_{\Zm}
\leq \frac{1}{T_{i+1}(\omega)^{1/n}}.
$$
Pour  montrer l'in\'egalit\'e de gauche,
consid\'erons des vecteurs entiers $w_i$ 
tels que
$$\|T_i(\omega)\omega \|_{\Zm}=
|T_i(\omega)\omega-w_i|.
$$
Comme 
$
T_i(\omega)w_{i+1}-T_{i+1}(\omega)w_i
$
est un vecteur entier non nul, on a
\begin{align*}
\frac{1}{T_i(\omega)T_{i+1}(\omega)}
& \leq 
\left|
\frac{T_i(\omega)w_{i+1}-T_{i+1}(\omega)w_i}{T_i(\omega)T_{i+1}(\omega)}
\right|
=\left|
\frac{w_{i+1}}{T_{i+1}(\omega)}-\frac{w_{i}}{T_{i}(\omega)}
\right|\\
& \leq
\left|
\frac{w_{i+1}}{T_{i+1}(\omega)}-\omega
\right|+
\left|
\omega-\frac{w_{i}}{T_{i}(\omega)}
\right|\\
& 
\leq 
\frac{\|T_{i+1}(\omega)\omega\|_{\Zm}}{T_{i+1}(\omega)}
+
\frac{\|T_{i}(\omega)\omega\|_{\Zm}}{T_{i}(\omega)}\\
&<
\frac{T_i(\omega)+T_{i+1}(\omega)}
{T_i(\omega)T_{i+1}(\omega)}
\|T_{i}(\omega)\omega\|_{\Zm},
\end{align*}
la Propri\'et\'e en d\'ecoule.

\section{Propri\'et\'es de transfert}
\subsection{ }
Rappelons que la suite $T_i(\omega)$ est la suite des
p\'eriodes introduite dans l'introduction.
On d\'efinit les ensembles
\begin{align*}
\Omega_n(\tau) &=
\left\{
\omega \in \Rm^n /\;\exists C>0, 
\forall T\in \Nm, \|T\omega\|_{\Zm}\geq 
CT^{-(1+\tau)/n} \right\},\\
\Omega^n(\tau) &=
\left\{
\omega \in \Rm^n /\;\exists C>0, 
\forall k\in \Zm^n-\{0\},\;
 \|\langle k,\omega\rangle\|_{\Zm}\geq 
C|k|^{-(1+\tau)n} \right\},\\
\Omega(\tau) &=
\left\{
\omega \in \Rm^n /\;\exists C>0, 
\forall i\in \Nm,\;
T_{i+1}(\omega)\leq CT_i(\omega)^{1+\tau}
\right\},\\
\tilde \Omega(\tau)&=
\left\{
\omega \in \Omega(\tau) / \;
\forall k\in \Zm^n-\{0\},\;
\langle k,\omega\rangle \not\in \Zm 
\right\}.
\end{align*}
Au vu  des th\'eor\`emes de Dirichlet de la section \ref{Dirichlet},
les ensembles $\Omega^n(\tau)$ et  $\Omega_n(\tau)$
sont vides pour $\tau<0$, il en va \'evidemment de m\^eme
des ensembles $\Omega(\tau)$ et
$\tilde \Omega(\tau)$.
On dit que les \'el\'ements de $\Omega^n(\tau)$
 pour $\tau \geq 0$ satisfont une propri\'et\'e diophantienne
lin\'eaire, et que les \'el\'ements de $\Omega_n(\tau)$
satisfont une propri\'et\'e  diophantienne simultan\'ee.
Dans cette note, on montre le

\subsection{ }
\textsc{Th\'eor\`eme de transfert : }Pour tout $\tau\geq 0$,
$$
\Omega_n\left(
\frac{\tau}{(n-1)\tau+n}\right)
\subset
\tilde \Omega\left(
\frac{\tau}{(n-1)\tau+n}\right)
\subset
\Omega^n(\tau)
\subset
\Omega_n(n\tau).
$$
On en d\'eduit par exemple 
que toutes les notions de vecteurs mal approchables
co\" \i ncident:
$$
\Omega_n(0)=\Omega^n(0)=\tilde\Omega(0).
$$ 
Lorsque $n=1$, ces nombres sont dits 
\textit{de type constant}.
On  d\'eduit aussi qu'une propri\'et\'e diophantienne
lin\'eaire implique toujours une propri\'et\'e diophantienne
simultan\'ee. L'apparten\-ance \`a $\Omega_n(\tau)$
n'implique cependant une propri\'et\'e diophantienne lin\'eaire que si 
$\tau<(n-1)^{-1}$. L'existence de ce seuil sera expliqu\'ee
dans la section \ref{rationnel}.
\subsection{ }\label{transfert}
Mentionnons que les inclusions suivantes, dues \`a Khintchine,
sont bien connues depuis longtemps
 (voir \cite{Lochak}, \cite{Schmidt}):
Pour tout $\tau\geq 0$,
\begin{equation*}
\Omega_n\left(
\frac{\tau}{(n-1)\tau+n}\right)
\subset
\Omega^n(\tau)
\subset
\Omega_n(n\tau).
\end{equation*}
%
%
%
%
%
%
%
%
%
%
%
%
%
\subsection{ }
La section suivante est consacr\'ee \`a quelques 
discussions sur le th\'eor\`eme de transfert.
On y r\'eduit  la preuve 
des deux nouvelles inclusions \`a celle 
la Proposition \ref{estimee}.
\section{Vecteurs rationnels et r\'esonances}\label{rationnel}

\subsection{ }
Le module de r\'esonance 
de $\omega$,
$R(\omega)$, est l'ensemble des vecteurs $k \in \Zm^n$
pour lesquels $\langle k,\omega\rangle\in\Zm$.
On appelle ordre de r\'esonance de $\omega$
le rang de ce sous groupe de $\Zm^n$,
c'est \`a dire le cardinal de ses bases.
Le vecteur $\omega$ est dit r\'esonant si son module de
r\'esonance n'est pas trivial, c'est \`a dire si son ordre de 
r\'esonance est non nul.
%
%
%
%
%
\subsection{ }
Le vecteur $\omega$ est dit rationnel si ses composantes 
sont rationnelles. Le vecteur $\omega$ est rationnel si et
seulement si toutes les orbites des $t_{\omega}$
sont p\'eriodiques.  Elles ont alors toutes la m\^eme
p\'eriode, qui est le plus petit d\'enominateur
commun des composantes de $\omega$.
Les vecteurs rationnels sont les vecteurs
maximalement r\'esonants, c'est \`a dire ceux dont l'ordre de
r\'esonance est $n$.
En effet, tout module $R$ de rang $n$ contient le module
$T\Zm^n$ pour une certain $T$. Les vecteurs 
ayant $R$ pour module de r\'esonance
sont donc n\'ecessairement $T$-p\'eriodiques.%
\subsection{ }
\textsc{Proposition : }\label{normal}
Soit $\omega$ un vecteur r\'esonant d'ordre $r$.
Il existe une  matrice $A\in Gl_n(\Zm)$
et une famille $d_i$, $1\leq i\leq r$, d'entiers
tels que $d_i$ divise $d_{i+1}$
et tels que 
$$
(d_1 e_1,\ldots, d_r e_r)
$$
forme une base  du module de r\'esonance de $A\omega$,
o\`u $(e_1,\ldots,e_n)$ est la base standard de $\Rm^n$.
En particulier,
$$
A\omega=(w/T,\omega')\in \Qm^{r}\times\Rm^{n-r},
$$
o\`u $\omega'$ est un vecteur non r\'esonant
et $w/T$ est un vecteur p\'eriodique (rationnel) de p\'eriode 
$T=d_r$.
Lorsque $\omega$ est p\'eriodique,
on lui associe donc $n$ invariants $d_i$,
o\`u  $d_n=T$ est la p\'eriode.
Il existe alors une  matrice   $A\in Gl_n(\Zm)$
telle que 
$$A\omega
=
(a_1/d_1, a_2/d_2, \ldots, a_n/d_n)
$$
avec des composantes $a_i/d_i$ irr\'eductibles.

\subsection{ }
\textsc{D\'emonstration : }
Soit $R$ le module de r\'esonance de $\omega$.
Consid\'erons une base $k_1,\ldots,k_r$ de $R$,
et associons lui la matrice $B$  \`a coefficients 
entiers 
de l'application lin\'eaire 
$\Rm^n\ni x\lto (\langle k_i,x\rangle)\in \Rm^r$.
Par un r\'esultat classique, voir \cite{Jac}, th\'eor\`eme 3.8,
il existe une matrice $A\in Gl_n(\Zm)$, une matrice 
$C\in Gl_r(\Zm)$ et une matrice diagonale 
$$
\Lambda=
\left[
\begin{matrix}
{d_1} & &  &0&0 &\cdots\\
& &\ddots & &&\\
0& & & d_r &0& \ldots\\
\end{matrix}
\right]
$$
telle que $B=C\Lambda A$,
o\`u les coefficients diagonaux $d_i$
sont des entiers tels que $d_i$ divide
$d_{i+1}$.
La propostion en d\'ecoule puisque
$$
\langle k,A \omega \rangle\in \Zm
\eq
\langle A^t k,\omega\rangle\in\Zm
\eq
A^t k\in \text{Im}(B^t)
=
\text{Im}(A^t\Lambda^t)
\eq 
k\in \text{Im}(\Lambda^t).
$$

\subsection{ }
\textsc{Remarque : }
Tout sous groupe  de $\Zm^n$ n'est pas un module
de r\'esonance, et les invariants d'un module de
r\'esonance satisfont certaines contraintes.
En particulier, il est facile de v\'erifier
que les invariants qui ne sont pas \'egaux \`a 
$1$ sont tous distincts:
$d_i=d_{i+1} \Rightarrow d_i=1$.

%
%
%
%
%
%
%
%
%
\subsection{ }
Lorsque  $n=1$, un rationnel $\omega=w/T$
(sous forme irr\'eductible) de grand
d\'enominateur $T$ 
a un comportement tr\`es proche d'un irrationnel 
en ce sens que les points $k\omega$, $k\in \Zm$
sont bien r\'epartis sur le cercle
(ce sont pr\'ecis\'ement les points de la forme
$l/T$, $l\in \Zm$).
En dimension sup\'erieure,
certains vecteurs p\'eriodiques de grande
p\'eriode se comportent essentiellement comme
des vecteurs non-r\'esonants, et d'autres 
se comportent plut\^ot comme des vecteurs 
r\'esonants (par exemple le  vecteur $(0,w/T)\in \Rm^2$).
En raison de la non compacit\'e
du groupe $Gl_n(\Zm)$ pour $n\geq 2$,
il n'est pas possible de quantifier 
la pr\'esence effective ou non de r\'esonances 
pour un vecteur p\'eriodique $\omega$ de grande p\'eriode
par des quantit\'es invariantes par l'action de 
$Gl_n(\Zm)$.
Par exemple, \'etant donn\'e un vecteur
$\omega\in  \Rm^2$,
l'algorithme des fractions continues
permet de trouver une suite (non born\'ee)
$A_n$ de matrices de $Gl_2(\Zm)$
et une suite de rationnels $x_n \in \Qm$
tels que $A_n(0,x_n)\lto \omega$
lorsque $n$ tend vers l'infini.

Il est donc utile d'introduire la quantit\'e
$$
e(\omega)=
\min
_{k\in R(\omega)-\{0\}}
|k|.
$$
Pour illustrer son  r\^ole, consid\'erons
un vecteur 
r\'eel $\omega$ et une suite d'approximations rationnelles
$w_n/T_n$ de $\omega$, o\`u la p\'eriode $T_n$ tend
vers l'infini. Il n'est pas difficile de voir que 
la limite $\omega$ est r\'esonante si et seulement si
la suite $e(w_n/T_n)$ est born\'ee.
Notons que, pour un vecteur rationnel $w/T$
de p\'eriode $T$, on a l'estimation
$$
1\leq e(w/T)\leq T^{1/n}.
$$
En effet, le th\'eor\`eme
\ref{lin}
donne l'existence d'un vecteur entier $k$ tel que
$|k|\leq T^{1/n}$ et 
$\| \langle k, w/T \rangle \|_{\Zm} <  1/T,$
mais alors $k\in R(w/T)$.
R\'eciproquement, pour la plupart des vecteurs rationnels $w/T$ de 
d\'enominateur $T$, $e(w/T)$ est au moins  de l'ordre de $T^{1/(n+1)}$.
Ceci est rendu plus pr\'ecis par la Propri\'et\'e \ref{e} ci-dessous.
Les vecteurs p\'eriodiques pour lesquels 
la valeur $e$ est grande sont ceux dont l'orbite
est constitu\'ee de points bien r\'epartis sur le tore.

\subsection{ }\label{e}
\textsc{Propri\'et\'e : }
Pour tout r\'eel strictement positif $\tau$,
la proportion de vecteurs entiers $w$ tels que 
$$
e(w/T)\leq T^{\frac{1-\tau}{n+1}}
$$
tend vers zero lorsque $T$ tend vers l'infini.

\textsc{D\'emonstration : }
Nous allons majorer le nombre de points entiers $w$
de $[0,T-1]^n$ tels qu'il existe un vecteur entier $k$
v\'erifiant  $|k|\leq A$
et $\langle k, w \rangle =Tl$, $l\in \Zm$.
Pour chaque couple $(k,l)\in \Zm^n\times\Zm$,
l'hyperplan  $\langle k, w \rangle =Tl$ contient au plus 
$T^{n-1}$ points entiers dans  $[0,T-1]^n$.
Par ailleurs, pour un vecteur entier 
$k$ donn\'e, 
le nombre d'entiers $l$ tels que 
l'hyperplan  $\langle k, w \rangle =Tl$
intersecte  $[0,T-1]^n$ est au plus 
$n(2|k|+1),$
puisque $|\langle k, w \rangle|\leq n|k||w|.$
Le nombre de couples $(k,l)\in \Zm^n\times\Zm$ 
v\'erifiant $|k|\leq A$ et
tels que l'hyperplan  d'\'equation $\langle k, w \rangle =Tl$
intersecte  $[0,T-1]^n$ est donc inf\'erieur \`a 
$n(2A+1)^{n+1}.$
Il y a donc au plus
$$
n(2A+1)^{n+1}T^{n-1}
$$
points entiers $w$  dans  $[0,T-1]^n$ tels que 
$e(w/T)\leq A.$
On termine alors la d\'emonstration en remarquant que
$$
n\left(
2 T^{\frac{1-\tau}{n+1}}+1
\right)^{n+1}
T^{n-1}
=o(T^n).
$$

\subsection{ }
En ne retenant d'un vecteur rationnel que sa p\'eriode $T$,
on perd donc  une information essentielle.
Par exemple, la suite des distances  $\|T(\omega,0)\|_{\Zm}$, $T\in \Zm$,
associ\'ee au vecteur $(\omega,0)\in \Rm^{n+1}$
est la m\^eme que la suite
 $\|T\omega\|_{\Zm}$
associ\'ee \`a $\omega$, alors que le vecteur $(\omega,0)$
est fortement r\'esonant et engendre des orbites
tr\`es mal r\'eparties dans $\Tm^{n+1}$.
De la m\^eme fa\c con, la suite des p\'eriodes
$T_i(\omega,0)$ est la m\^eme que celle de $\omega$.
Il est  remarquable  qu'une simple connaissance 
de la suite  $\|T\omega\|_{\Zm}$ puisse suffire
\`a d\'eterminer le caract\`ere non r\'esonant de $\omega$,
voir \ref{nonres}.
Il n'est pas possible de tirer ce caract\`ere non r\'esonant d'une 
connaissance de la suite des p\'eriodes.
Pourtant, si l'on suppose \`a priori le vecteur $\omega$  
non r\'esonant,
on montre que certaines estim\'ees sur la croissance des p\'eriodes 
permettent de d\'eduire des propri\'et\'es diophantiennes lin\'eaires
v\'erifi\'ees par $\omega$. Cette remarque,
\`a ma connaissance nouvelle,
fait l'objet de la proposition \ref{estimee},
qui sera d\'emontr\'ee dans la section suivante.

%
%
%
%
%
%
%
%

\subsection{ }\label{nonres}
\textsc{Propri\'et\'e : }Pour tout $\tau<(n-1)^{-1}$,
les vecteurs de  $\Omega_n(\tau)$
sont non r\'esonants.\\
Ce r\'esultat est bien sur impliqu\'e par les inclusions 
de Khintchine \ref{transfert}. Nous allons 
cependant en donner une preuve 
\'el\'ementaire.
Consid\'erons un vecteur r\'esonant $\omega.$
D'apr\`es \ref{normal}, il existe une matrice 
$A \in SL_n(\Zm)$ telle que 
$A\omega=(\omega',w/p) \in \Rm^{n-1} \times \Qm $, $(w,p)\in\Zm^2$.
Appliquons le th\'eor\`eme de Dirichlet \`a
$p\omega'$.
Il existe une infinit\'e d'entiers $T$
tels que
$\|Tp\omega'\|_{\Zm} <T^{-1/{n-1}}$, et donc tels que 
$$
\|T(p\omega',w)\|_{\Zm}
< \frac{1}{T^{1/(n-1)}}
\eq
\|TpA\omega\|_{\Zm}
<\frac{p^{1/(n-1)}}{(pT)^{1/(n-1)}}.
$$
Il y a donc une infinit\'e d'entiers $T'=pT$ et une constante $C$ 
tels que 
$$
\|T'\omega\|_{\Zm}
\leq\frac{C}{(T')^{1/(n-1)}}.
$$
Ceci exclut que $\omega$ puisse appartenir
\`a l'ensemble $\Omega_n(\tau)$
lorsque
$$
\frac{1+\tau}{n}<\frac{1}{n}
\eq
\tau <
\frac{1}{n-1}.
$$

\subsection{ }\label{estimee}
\textsc{Proposition :}
Si $\omega$ est non r\'esonant  et si il existe 
$\tau \in \left[ 0,\frac{1}{n-1}\right[$ et $C>0$
tels que
$
T_{i+1}(\omega)\leq C T_i(\omega)^{1+\tau}
$
alors il y a une constante $D>0$ telle que, pour tout 
$k\in \Zm^n$, on a 
$$
\|\langle k,\omega \rangle \|_{\Zm}\geq D|k|^
{-\frac{1+\mu}{n}},\;
\text{   o\`u   } \,\;
\mu= \frac{n\tau}{n-(n-1)(1+\tau)}.
$$

\subsection{ }
\textsc{D\'emonstration du th\'eor\`eme de transfert :}
Nous ne prouverons que les deux premi\`ere inclusions.
Nous renvoyons \`a \cite{Schmidt} pour la preuve de la troisi\`eme
inclusion.
La deuxi\`eme inclusion d\'ecoule directement
de la propri\'et\'e \ref{estimee}.
La premi\`ere inclusion d\'ecoule de la propri\'et\'e
\ref{nonres}
et de la remarque suivante :
Si $\omega\in \Omega_n(\tau')$,
on obtient en appliquant  th\'eor\`eme de Dirichlet
avec $Q=T_{i+1}(\omega)$
que 
$$
CT_i(\omega)^{-(1+\tau')/n}\leq \|T_i(\omega)\omega\|_{\Zm}
\leq T_{i+1}(\omega)^{-1/n},
$$
et donc que $\omega\in \Omega(\tau')$
 (voir \cite{Lochak}). 

\section{D\'emonstration de la proposition \ref{estimee}}
Le vecteur $\omega$ est fix\'e une fois pour toutes,
on notera donc $T_i$ les p\'eriodes $T_i(\omega)$.
On consid\`ere une  suite  de meilleures approximations $\omega_i$,
c'est \`a dire une suite de vecteurs rationnels 
telle que $\|T_i (\omega_i-\omega)\|=\|T_i \omega\|_{\Zm}$
et 
$T_i\omega_i=w_i\in \Zm^n.$
Soit $k\in \Zm^n$, comme $\omega$ est non r\'esonant,
les produits scalaires $\langle k,\omega_i \rangle$
ne peuvent \^etre entiers que pour un nombre fini 
de valeurs de $i$, et on peut d\'efinir
$$
i(k)=1+\max 
\big\{ i \text{ tel que } \langle k,\omega _i\rangle\in \Zm\big\}.
$$
\subsection{ }
\'Etablissons pour commencer les in\'egalit\'es 
$$
\|k\|\geq
\frac{T_{i(k)-1}}{2\sqrt{n}\,T_{i(k)}^{1-1/n}}
\geq
\frac{1}{2\sqrt{n}\, C^{\frac{1}{1+\tau}}}
 T_{i(k)}^{\frac{n+(1-n)(1+\tau)}{n(1+\tau)}}.
$$
On a 
\begin{align*}
 \langle k,\omega _{i(k)-1}\rangle &\in  \Zm,\\
\intertext{et}
\|\langle k,\omega _{i(k)}\rangle\|_{\Zm} &\geq 1/T_{i(k)}\\
\intertext{donc}
|\langle k,\omega _{i(k)}-\omega _{i(k)-1}\rangle| 
&\geq 1/T_{i(k)}.
\end{align*}
On en d\'eduit que 
$$
\|k\|\geq \frac{1}{2T_{i(k)}\|\omega_{i(k)-1}-\omega\|}.
$$
En appliquant 
le Th\'eor\`eme de Dirichlet
avec $Q=T_{i+1}$, on obtient
$$
\|\omega_{i(k)-1}-\omega\|
\leq \sqrt{n}|\omega_{i(k)-1}-\omega|
\leq  \frac{ \sqrt{n}}
{T_{i(k)-1}T_{i(k)}^{1/n}},
$$
ce qui donne la premi\`ere \'egalit\'e.
La seconde in\'equalit\'e est alors une cons\'equence directe de 
l'hypoth\`ese.

\subsection{ }
Un calcul simple donne alors que 
$$
T_{i(k)}\leq \big(2\sqrt{n}\, C^{\frac{1}{1+\tau}}
\|k\|\big)^{n(1+\mu)},
$$
o\`u $\mu$ est donn\'e dans l'\'enonc\'e.
D\'efinissons l'indice
$$
j(k)=\min \big\{ j\geq i(k) \text{ tels que }
T_{j+1}\geq  \big(2\sqrt{n}\, C^{\frac{1}{1+\tau}}
\|k\|\big)^{n(1+\mu)}\big\}.
$$
D'un cot\'e, l'in\'egalit\'e
$$
T_{j(k)+1}\geq  \big(2\sqrt{n}\, C^{\frac{1}{1+\tau}}
\|k\|\big)^{n(1+\mu)}
$$
implique, puisque  $T_{j(k)+1}\leq CT_{j(k)}^{1+\tau},$ que
$$
T_{j(k)}\geq 
\frac{ \big(2\sqrt{n}\, C^{\frac{1}{1+\tau}})^{n\frac{1+\mu}{1+\tau}}}
{C^{\frac{1}{1+\tau}}}\,
\|k\|^{n\frac{1+\mu}{1+\tau}}.
$$
D'un autre cot\'e, on a la majoration
$$
T_{j(k)}\leq  \big(2\sqrt{n}\, C^{\frac{1}{1+\tau}}
\|k\|\big)^{n(1+\mu)}.
$$
On a donc, comme $j(k)\geq i(k)$,
\begin{equation}\label{min}
\|\langle k, \omega_{j(k)}\rangle \|_{\Zm}
\geq \frac{1}{T_{j(k)}}
\geq \big(2\sqrt{n}\, C^{\frac{1}{1+\tau}}
\|k\|\big)^{-n(1+\mu)}.
\end{equation}
En utilisant encore le th\'eor\`eme de Dirichlet, on obtient
$$
|\omega_{j(k)}-\omega|\leq \frac{1}{T_{j(k)}^{1+1/n}}
$$
et donc 
\begin{equation}\label{delta}
|\langle k, \omega_{j(k)}-\omega\rangle|
\leq\frac{\|k\|}{T_{j(k)}^{1+1/n}}
\leq  \frac{1}{T_{j(k)}}\tilde C 
\|k\|^{\frac{\tau -\mu}{1+\tau}},
\end{equation}
o\`u $\tilde C$
 ne d\'epend pas de $k$.
Il est facile de voir que $\tau<\mu$,
et donc, par (\ref{min})  et (\ref{delta})  que 
$$
\|\langle k,\omega\rangle \|_{\Zm}\geq
\frac{1}{T_{j(k)}}(1-\epsilon(\|k\|))
\geq \big(2\sqrt{n}\, C^{\frac{1}{1+\tau}}
\|k\|\big)^{-n(1+\mu)}(1-\epsilon(\|k\|)),
$$
o\`u $\epsilon(x)$ est une fonction qui tend vers $0$ lorsque $x$
tend vers l'infini.
La proposition annonc\'ee en d\'ecoule.


\begin{thebibliography}{99} 

\bibitem{Jac} Jacobson N.: \textit{Basic Algebra I,} 
second edition,
W. H. Freeman and Company, New York, (1985).


\bibitem{Lochak} Lochak P. :
Canonical perturbation theory via simultaneous approximation,
Russ. Math. Surveys \textbf{47} (1992) 57-133.

\bibitem{Schmidt} Schmidt W. M. : \textit{Diophantine Approximation},
Lecture Notes in Math. \textbf{785} (1980).
\end{thebibliography}
\end{document}